\input amstex
\documentstyle{amsppt}
\magnification=1200
\input scrload
\define \C#1{{\scr #1}}
\define \my{\bf}
\define\po{\parindent 0pt}
\define \mytitle{\po \m \my}
\define \p{{\po \it \m Proof. }}
\define \m{\medskip}
\nologo
\define \ph#1{\phantom {#1}}
\define \T#1{\widetilde {#1}}
\define \Ker{\operatorname{Ker}}
\define \IM{\operatorname {Im}}
\define \Fix{\operatorname {Fix}}
\define \cat{\operatorname{cat}}
\define \inter{\operatorname{int}}
\define \cl{\operatorname{cl}}
\define \Crit{\operatorname{Crit}}
\define \crit{\operatorname{crit}}

\define \pt{\operatorname{pt}}
\define \Arn{\operatorname{Arn}}
\define \Rest{\operatorname{Rest}}

\define \mul{{\raise0.2ex \hbox{$\scriptstyle \bullet$}}}
\topmatter
\title{On analytical applications of stable homotopy (the Arnold  
conjecture, critical points) }\endtitle
\rightheadtext{On analytical applications of stable homotopy}
\author{Yuli B. Rudyak}\endauthor
\address{Mathematisches Institut,
 Universit\"at Heidelberg,
Im Neuenheimer Feld 288,
D-69120 Heidelberg 1,
Germany}
\endaddress
\email{july\@mathi.uni-heidelberg.de}\endemail
\thanks{The author was partially supported by Deutsche  
Forschungsgemeinschaft, the Research Group ``Topologie und  
nichtkommutative Geometrie''}\endthanks
\abstract  We prove the Arnold conjecture for closed symplectic  
manifolds with $\pi_2(M)=0$ and $\cat M=\dim M$. Furthermore, we  
prove an analog of the Lusternik--Schnirelmann theorem for functions  
with ``generalized hyperbolicity'' property. \endabstract
\endtopmatter
\head {\bf Introduction} \endhead
Here we show that the technique developed in [R98] can be applied  
to the Arnold conjecture and to estimation of the number of critical  
points. For convenience of the reader, this paper is written  
independently of [R98].
\m Given a smooth ($=C^{\infty}$) manifold $M$ and a smooth  
function $f: M \to \Bbb R$, we denote by $\crit f$ the number of  
critical points of $f$ and set $\Crit M=\min \{\crit f\}$ where $f$  
runs over all smooth functions $M \to \Bbb R$.
\m The Arnold conjecture~[Ar89, Appendix 9] is a well-known problem  
in Hamiltonian dynamics. We recall the formulation. Let  
$(M,\omega)$ be a closed symplectic manifold, and let $\phi: M \to  
M$ be a Hamiltonian symplectomorphism (see [HZ94], [MS95] for the  
definition). Furthermore, let $\Fix \phi$ denote the number of fixed  
points of $\phi$. Finally, let
$$
\Arn(M,\omega):=\min_\phi{\Fix\phi}
$$
where $\phi$ runs over all Hamiltonian symplectomorphisms $M \to M$.
The Arnold conjecture claims that $\Arn(M,\omega)\geq \Crit (M)$.
\par It is well known and easy to see that $\Arn(M,\omega) \leq  
\Crit M$. Thus, in fact, the Arnold conjecture claims the equality  
$\Arn (M, \omega)=\Crit M$.
\m Let $\cat X$ denote the Lusternik--Schnirelmann category of a  
topological space $X$ (normalized, i.e., $\cat X=0$ for $X$  
contractible).
\m Given a symplectic manifold $(M^{2n}, \omega)$, we define the  
homomorphisms
$$
\split
I_{\omega}:\pi_2(M) \to \Bbb Q,&\quad I_{\omega}(x)=\langle \omega,  
h(x) \rangle \\
I_c:\pi_2(M) \to \Bbb Z,&\quad I_c(x)=\langle c, h(x) \rangle
\endsplit
$$
where $h: \pi_2(M) \to H_2(M)$ is the Hurewicz homomorphism,  
$c=c_1(\tau M)$ is the first Chern class of $M$ and $\langle  
-,-\rangle $ is the Kronecker pairing.
\proclaim{Theorem A {\rm (see 3.6)}} Let $(M^n,\omega), n=\dim M$  
be a closed connected symplectic manifold such that  
$I_{\omega}=0=I_c$ $($e.g., $\pi_2(M)=0)$ and $\cat M=n$. Then  
$\Arn(M,\omega)\geq \Crit M$, i.e., the Arnold conjecture holds for  
$M$.
\endproclaim
\m It is well known that $\Crit M \geq 1+\cl M$ for every closed  
manifold $M$, where $\cl$ denotes the cup-length, i.e., the length  
of the longest non-trivial cup-product in $\T H^*(M)$. So, one has  
the following weaker version of the Arnold conjecture:
$$
\Arn (M,\omega) \geq 1+\cl (M),
$$
and most known results deal with this weak conjecture, see [CZ83],  
[S85], [H88], [F89-1], [F89-2], [LO96]. (Certainly, there are lucky  
cases when $\Crit M =1 +\cl M$, e.g. $M=T^{2n}$, cf. [CZ83].) For  
example, Floer~[F89-1],~[F89-2] proved that $\Arn (M) \geq 1+\cl M$  
provided $I_{\omega}=0=I_c$, cf. also Hofer~[H88]. So, my  
contribution is the elimination of the clearance between $\Crit M$  
and $1+\cl M$. (It is easy to see that there are manifolds $M$ as in  
Theorem A with $\Crit M > 1+\cl M$, see 3.7 below.)
\m Actually, we prove that $\Arn(M,\omega)\geq 1+\cat M$ and use a  
result of Takens~[T68] which implies that $\Crit M=1+\cat M$  
provided $\cat M =\dim M$.
\m After submission of the paper the author and John Oprea proved  
that $\cat M= \dim M$ for every closed symplectic manifold  
$(M,\omega)$ with $I_{\omega}=0$, see [RO97]. So, the condition  
$\cat M =n$ in Theorem A can be omitted.
\m Passing to critical points, we prove the following theorem.
\proclaim{Theorem B {\rm (see 4.5)}} Let $M$ be a closed orientable  
manifold, and let $g: M \times \Bbb R^{p+q} \to \Bbb R$ be a  
$C^2$-function with the following properties:
\roster
\item There exist disks $D_+\i \Bbb R^p$ and $D_-\i \Bbb R^q$  
centered in origin such that $\inter (M \times D_+ \times D_-)$  
contains all critical points of $g$;
\item $\nabla g(x)$ points inward on $M \times \partial D_+ \times  
\inter D_-$ and outward on $M \times \inter D_+ \times \partial  
D_-$.
\endroster
Then $\crit g\geq 1+r(M)$.
\par In  particular, if $M$ is aspherical then $\crit g \geq 1+\cat  
M$. \footnote"$^\dagger$" {Recall that a connected topological  
space is called aspherical if $\pi_i(X)=0$ for every $i>1$.}
\endproclaim
\m Notice that functions $g$ as in Theorem B are related to the  
Conley index theory, see [C76]. I  remark that Cornea~[Co98] have  
also estimated the number of crirical points of functions as in  
Theorem B.
\m We reserve the term ``map'' for continuous functions of  
topological spaces, and we call a map {\it inessential} if it is  
homotopic to a constant map. The disjoint union of spaces $X$ and  
$Y$ is denoted by $X\sqcup Y$. Furthermore, $X^+$ denotes the  
disjoint union of $X$ and a point, and $X^+$ is usually considered  
as a pointed space where the base point is the added point.
\m We follow Switzer~[Sw75] in the definition of $CW$-complexes. A  
{\it $CW$-space} is defined to be a space which is homeomorphic to a  
$CW$-complex.
\m Given a pointed $CW$-complex $X$, we denote by  
$\Sigma^{\infty}X$ the spectrum $E=\{E_n\}$ where $E_n=S^nX$ for  
every $n\geq 0$ and $E_n=\pt$ for $n<0$; here $S^nX$ is the $n$-fold  
reduced suspension over $X$. Clearly, $\Sigma^{\infty}$ is a  
functor from pointed $CW$-complexes to spectra.
\m Given any (bad) space $X$, the cohomology group $H^n(X;\pi)$ is  
always defined to be the group $[X,K(\pi,n)]$ where $[-,-]$ denotes  
the set of homotopy classes of maps.
\m ``Smooth'' always means ``$C^{\infty}$''.
\par ``Fibration'' always means a Hurewicz fibration.
\par ``Connected'' always means path connected.
\par The sign `` $\simeq $'' denotes homotopy of maps (morphisms)  
or homotopy equivalence of spaces (spectra).
\head {\bf \S 1. Preliminaries on the Lusternik--Schnirelmann  
category} \endhead
{\mytitle 1.1. Definition.} (a) Given a subspace $A$ of a  
topological space $X$, we define $\cat_XA$ to be the minimal number  
$k$ such that $A= U_1\cup \cdots \cup U_{k+1}$ where each $U_i$ is  
open in $A$ and contractible in $X$. We also define $\cat_XA=-1$ if  
$A=\emptyset$.
\par (b) Given a map $f: X \to Y$, we define $\cat f$ to be the  
minimal number $k$ such that $X= U_1\cup \cdots \cup U_{k+1}$ where  
each $U_i$ is open in $X$ and $f|U_i$ is inessential for every $i$.
\par (c) We define the {\it Lusternik--Schnirelmann category} $\cat  
X$ of a space $X$ by setting $\cat X:= \cat_XX=\cat 1_X$.
\par Clearly, $\cat f \leq \min\{\cat X, \cat Y\}$.
\m The basic information concerning the Lusternik--Schnirelmann  
category can be found in [Fox41], [J78], [Sv66].
\m Let $X$ be a connected space. Take a point $x_0\in X$, set
$$
PX=P(X,x_0)=\{\omega\in X^I\bigm|\omega (0)=x_0\}
$$
and consider the fibration $p: PX \to X,\, p(\omega)=\omega(1)$  
with the fiber $\Omega X$.
\m Given a natural number $k$, we use the short notation
$$
p_k: P_k(X) \to X.
\tag{1.2}
$$
for the map
$$
\underbrace{p_X*_X \cdots *_X p_X}_{k\text{ times}}:  
\underbrace{PX*_X\cdots *_X PX}_{k\text{ times}} @>\ph{OOOO}>> X
$$
where $*_X$ denotes the fiberwise join over $X$, see e.g. [J78].
In  particular, $P_1(X)=PX$.
\proclaim{1.3. Proposition} For every connected compact metric  
space $X$ and every natural number $k$ the following hold: \par
{\rm (i)} $p_k: P_k(X) \to X$ is a fibration; \par
{\rm (ii)} $\cat P_k(X) < k$; \par
{\rm (iii)} The homotopy fiber of the fibration $p_k: P_k(X) \to X$  
is the $k$-fold join $(\Omega X)^{*k}$; \par
{\rm (iv)} If $\cat X=k$ and $X$ is $(q-1)$-connected then $p_k:  
P_k(X) \to X$ is a $(kq-2)$-equivalence; \par
{\rm (v)} If $X$ has the homotopy type of a $CW$-space then  
$P_k(X)$ does. \par
\endproclaim
\p (i) This holds since a fiberwise join of fibrations is a  
fibration, see [CP86]. \par
(ii) It is easy to see that $\cat(E_1*_XE_2)\leq \cat E_1 +\cat  
E_2+1$ for every two maps $f_1:E_2 \to X$ and $f_2: E_2 \to X$. Now  
the result follows since $\cat P_1(X)=0$. \par
(iii) This holds since the homotopy fiber of $p_1$ is $\Omega X$. \par 
(iv) Recall that $A*B$ is $(a+b+2)$-connected if $A$ is  
$a$-connected and $B$ is $b$-connected. Now, $\Omega X$ is  
$(q-2)$-connected, and so the fiber $(\Omega X)^{*k}$ of $p_k$ is  
$(kq-2)$-connected. \par
(v) It is a well-known result of Milnor [M59] that $\Omega X$ has  
the homotopy type of a $CW$-space. Hence, the space $(\Omega  
X)^{*k}$ has it. Finally, the total space of any fibration has the  
homotopy type of a $CW$-space provided both the base and the fiber  
do, see e.g. [FP90, 5.4.2].
\qed
\proclaim{1.4. Theorem {\rm ([Sv66, Theorems 3 and 19$'$])}} Let  
$f: X \to Y$ be a map of connected compact metric spaces. Then    
$\cat f<k$ iff there is a map $g: X\to P_k(Y)$ such that $p_kg=f$.
\qed
\endproclaim
\head {\S 2. An invariant $r(X)$} \endhead
Consider the Puppe sequence
$$
P_m(X) @>p_m>> X @>j_m>> C_m(X):=C(p_m)
$$
where $p_m: P_m(X) \to X$ is the fibration (1.2) and $C(p_m)$ is  
the cone of $p_m$.
{\mytitle 2.1. Definition.} Given a connected space $X$, we set
$$
r(X):=\sup\{m|j_m\text{ is stably  essential}\}.
$$
(Recall that a map $A \to B$ is called stably essential if it is  
not stably homotopic to a constant map.)
\proclaim{2.2. Proposition} {\rm (i)} $r(X) \leq \cat X$ for every  
connected compact metric space $X$. \par
{\rm (ii)} Let $X$ be a connected $CW$-space, let $E$ be a ring  
spectrum, and let $u_i\in \T E^*(X)), i=1, \ldots, n$ be elements  
such that $u_1 \cdots u_n \neq 0$. Then $r(X)\geq n$. In  other  
words, $r(X)\geq \cl_E(X)$ for every ring spectrum $E$.
\endproclaim
\m It makes sense to remark that $r(X)=\cat X$ iff $X$ possesses a  
detecting element, as defined in [R96].
\p (i) This follows from 1.4.\par
(ii) Because of 1.3(v), without loss of generality we can and shall  
assume that $C_n(X)$ is a $CW$-space. We set $u=u_1 \cdots u_n\in  
\T E^d(X)$. Because of the cup-length estimation of the  
Lusternik--Schnirelmann category, and by 1.3(ii), we have  
$p_n^*(u)=0$. Hence, there is a homotopy commutative diagram
$$
\CD
\Sigma^{\infty}X @>\Sigma^{\infty}j_n>> \Sigma^{\infty}C_n(X)\\
@VuVV @VVV\\
\Sigma^dE @= \,\Sigma^dE.
\endCD
$$
Now, if $r(X)<n$ then $\Sigma^{\infty}j_n$ is inessential, and so  
$u=0$. This is a contradiction.
\qed
\proclaim{2.3. Lemma} Let $f: X \to Y$ be a map of compact metric  
spaces with $Y$ connected, and let $j_m: Y \to C_m(Y)$ be as in  
$2.1$. If the map $j_mf$ is essential then $\cat f \geq m$. In   
particular, if $r(Y)=r$ and the map
$$
f^{\#}: [Y,C_r(Y)] \to [X,C_r(Y)]
$$
is injective then $\cat f\geq r$.
\endproclaim
\p Consider the diagram
$$
\CD
 @. X @. \\
@. @VfVV @.\\
P_m(Y) @>p_m>> Y @>j_m>> \,C_m(Y).
\endCD
$$
If $\cat f < m$ then, by 1.4, there is a map $g: X \to P_m(Y)$ with  
$f=p_mg$, and hence $j_m f$ is inessential. This is a  
contradiction.
\qed
\proclaim{2.4. Theorem} Let $M^n$ be a closed oriented connected  
$n$-dimensional PL manifold such that $\cat M=\dim M\geq 4$. Then  
$r(M)=\cat M$.
\endproclaim
\p Let $MSPL_*(-)$ denote the oriented PL bordism theory. By the  
definition of $r(X)$, it suffices to prove that $(j_n)_*:MSPL_*(M)  
\to MSPL_*(C_n(M))$ is a non-zero homomorphism. Hence, it suffices  
to prove that $(p_n)_*: MSPL_n(P_n(M)) \to MSPL_n(M)$ is not an  
epimorphism. Clearly, this will be proved if we prove that $[1_M]\in  
MSPL_n(M)$ does not belong to $\IM (p_n)_*$. \par
Suppose the contrary. Then there is a map $F: W \to M$ with the  
following properties:
\roster
\item $W$ is a compact $(n+1)$-dimensional oriented  PL manifold  
with $\partial W=M \sqcup  V$;
\item $F|M=1_M$, $F|V: V \to M$ lifts to $P_n(M)$ with respect to  
the map $p_n: P_n(M) \to M$.
\endroster
Without loss of generality we can assume that $W$ is connected.
\par Suppose for a moment that $\pi_1(W,M)$ is a one-point set  
(i.e., the pair $(W,M)$ is simply connected). Then $(W,M)$ has the  
handle presentation without handles of indices $\leq 1$, see [St68,  
8.3.3, Theorem A]. By duality, the pair $(W,V)$ has the handle  
presentation without handles of indices $\geq n$. In other words, $W  
\simeq V\cup e_1\cup \cdots\cup e_s$ where $e_1, \ldots, e_s$ are  
cells attached step by step and such that $\dim e_i\leq n-1$ for  
every $i=1, \ldots, s$. However, the fibration $p_n: P_n(M) \to M$  
is $n-2$ connected. Thus, $F: W \to M$ can be lifted to $P_n(M)$. In  
 particular, $p_n$ has a section. But this contradicts 1.4. \par
So, it remains to prove that, for every membrane $(W,F)$, we can  
always find a membrane $(U,G)$ with $\pi_1(U,G)=*$ and $G|\partial  
U=F|\partial W$. Here $\partial U = \partial W =M \sqcup V$ and $G:  
U \to M$. We start with an arbitrary connected membrane $(W,F)$.  
Consider a PL embedding $i: S^1 \to \inter W$. Then the normal  
bundle $\nu$ of this embedding is trivial. Indeed, $w_1(\nu)=0$  
because $W$ is orientable.
\par Since $M$ is a retract of $W$, there is a commutative diagram
$$
\CD
0 @>>> \pi_1(M) @>>>  \pi_1(W) @>>>\pi_1(W,M) @>>> 0 \\
@.@| @VVF_*V @.@.\\
@. \pi_1(M) @= \pi_1(M) @.@.
\endCD
$$
where the top line is the homotopy exact sequence of the pair  
$(W,M)$. Clearly, if $F_*$ is monic then $\pi_1(W,M)=*$.
\par Let $\pi_1(W)$ be generated by elements $a_1, \ldots, a_k$. We  
set $g_i:=F_*(a_i)a_i^{-1}\in \pi_1(W)$ where we regard $\pi_1(M)$  
as the subgroup of $\pi_1(W)$. Then $\Ker F_*$ is the smallest  
normal subgroup of $\pi_1(W)$ contained $g_1, \ldots, g_k$. Now we  
realize $g_1, \ldots, g_k$  by PL
embeddings $S^1 \to \inter W$ and perform the surgeries of $(W,F)$  
with respect to these embeddings, see [W70]. The result of the  
surgery establishes us a desired membrane.
\qed
\proclaim{2.5. Corollary} Let $M$ be as in $2.4$, let $X$ be a  
compact metric space, and let $f: X \to M$ be a map such that $f^*:  
H^n(M;\pi_n(C_n(M))) \to H^n(X;\pi_n(C_n(M)))$ is a monomorphism.  
Then $\cat f \geq \cat M$. \endproclaim
Notice that, in fact, $\cat f =\cat M$ since $\cat f \leq \cat M $  
for general reasons.
\p We set $\pi=\pi_n(C_n(M))$. It is easy to see that $C_n(M)$ is  
simply connected. Hence, by 1.3(iv) and the Hurewicz theorem,  
$C_n(M)$ is $(n-1)$-connected. Thus, $[M, C_n(M)]= H^n(M;\pi)$. Let  
$\iota: C_n(M) \to K(\pi,n)$ denote the fundamental class. Then  
$f^*$ can be decomposed as
$$
f^*: H^n(M;\pi) = [M, C_n(M)] @>f^{\#}>> [X, C_n(M)] @>\iota_*>>  
[X,K(\pi,n)]=H^n(X;\pi).
$$
Since $f^*$ is a monomorphism, we conclude that $f^{\#}$ is. Thus,  
by 2.3 and 2.4,
$$
\cat f \geq r(M)=\cat M.
\qed
$$
\head {\S 3. The invariant $r(M)$ and the Arnold conjecture} \endhead
\m Recall (see the introduction) that the Arnold conjecture claims  
that $\Arn(M,\omega)\geq \Crit M$ for every closed symplectic  
manifold $(M,\omega)$.
{\mytitle 3.1. Recollection.} A {\it flow} on a topological space  
$X$ is a family $\Phi=\{\varphi_t\}, t\in \Bbb R$ where each  
$\varphi_t: X \to X$ is a self-homeomorphism and  
$\varphi_s\varphi_t=\varphi_{s+t}$ for every $s,t\in \Bbb R$ (notice  
that this implies $\varphi_0=1_X$).
\par A flow is called {\it continuous} if the function $X \times  
\Bbb R \to X, (x,t)\mapsto \varphi_t(x)$ is continuous.
\par A point $x\in X$ is called a {\it rest point} of $\Phi$ if  
$\varphi_t(x)=x$ for every $t\in \Bbb R$. We denote by $\Rest \Phi$  
the number of rest points of $\Phi$.
\par A continuous flow $\Phi=\{\varphi_t\}$ is called {\it  
gradient-like} if there exists a continuous (Lyapunov) function $F:  
X \to \Bbb R$ with the following property: for every $x\in X$ we  
have $F(\varphi_t(x))< F(\varphi_s(x))$ whenever $t>s$ and $x$ is  
not a rest point of $\Phi$.
{\mytitle 3.2. Definition {\rm (cf. [H88], [MS95])}.} Let $X$ be a  
topological space. We define an {\it index function} on $X$ to be  
any function $\nu: 2^X \to \Bbb N\cup\{0\}$ with the following  
properties:
\roster \item (monotonicity) If $A\i B \i X$ then $\nu (A)\leq \nu (B)$;
\item (continuity) For every $A\i X$ there exists an open  
neighbourhood $U$ of $A$ such that $\nu (A) =\nu (U)$;
\item (subadditivity) $\nu(A\cup B) \leq \nu (A) + \nu (B)$;
\item (invariance) If $\{\varphi_t\}, t\in \Bbb R$ is a continuous  
flow on $X$ then $\nu (\varphi_t(A))=\nu(A)$ for every $A\i X$ and  
$t\in \Bbb R$;
\item (normalization) $\nu(\emptyset)=0$. Furthermore, if $A\neq  
\emptyset$ is a finite set which is contained in a connected  
component of $X$ then $\nu(A)=1$.
\endroster
\proclaim{3.3. Theorem} Let $\Phi$ be a gradient-like flow on a  
compact metric space $X$. Then
$$
\Rest \Phi \geq \nu(X)
$$
for every index function $\nu$ on $X$.
\endproclaim
\p The proof follows the ideas of Lusternik--Schnirelmann. For $X$  
connected see [H88], [MS95, p.346 ff]. Furthermore, if $X=\sqcup  
X_i$ with $X_i$ connected then
$$
\Rest \Phi=\sum \Rest(\Phi|X_i)\geq \sum \nu(X_i)\geq \nu(X).
\qed
$$
\proclaim{3.4. Corollary} Let $\Phi$ be a gradient-like flow on a  
compact metric space $X$, let $Y$ be a Hausdorff space which admits  
a covering $\{U_{\alpha}\}$ such that each $U_{\alpha}$ is open and  
contractible in $Y$, and let $f: X \to Y$ be an arbitrary map.  Then  

$$
\Rest \Phi \geq 1+\cat f.
$$
\endproclaim
\p Given a subspace $A$ of $X$, we define $\nu(A)$ to be the  
minimal number $m$ such that $A\i U_1\cup \cdots \cup U_m$ where  
each $U_i$ is open in $X$ and $f|U_i$ is inessential. It is easy to  
see that $\nu$ is an index function on $X$ (normalization follows  
from the properties of $Y$). But $\nu(X)=1+\cat f$, and so, by 3.3,  
we conclude that $\Rest \Phi \geq 1+\cat f$.
\qed
\proclaim{3.5. Theorem} Let $(M,\omega)$ be a closed connected  
symplectic manifold with $I_{\omega}=0=I_c$, and let $\phi: M \to M$  
be a Hamiltonian symplectomorphism. Then there exists a map $f: X  
\to M$ with the following properties: \par
{\rm (i)} $X$ is a compact metric space;\par
{\rm (ii)} $X$ possesses a gradient-like flow $\Phi$ such that  
$\Rest \Phi = \Fix \phi$; \par
{\rm (iii)} The homomorphism $f^*: H^n(M;G) \to H^n(X;G)$ is a  
monomorphism for every coefficient group $G$.
\endproclaim
\p This can be proved following [F89-2, Theorem 7]. (Note that the  
formulation of this theorem contains a misprint: there is typed  
$z^*[P] =0$, while it must be typed $z^*[P] \neq 0$. Furthermore,  
the reference [CE] in the proof must be replaced by [F7].) In fact,  
Floer denoted by $z : \C S \to P$ what we denote by $f: X \to M$,  
and he showed that the homomorphism $z^*: \Bbb Z=H^n(P) \to H^n(\C  
S)$ is monic. He did it for $\Bbb Z$-coefficients, but the proof for  
arbitrary $G$ is similar.
\par Also, cf. [H88] and [HZ94, Ch. 6].
\par In fact, Floer considered Alexander--Spanier cohomology, but  
for compact metric spaces it coincides with $H^*(-)$. In greater  
detail, you can find in [Sp66] an isomorphism between  
Alexander--Spanier and \v Cech cohomology and in [Hu61] an  
isomorphism between \v Cech cohomology and $H^*(-)$.
\qed
\m Recall that every smooth manifold turns out to be a PL manifold  
in a canonical way, see e.g~[Mu66].
\proclaim{3.6. Theorem} Let $(M,\omega)$ be a closed connected  
symplectic manifold with $I_{\omega}=0=I_c$ and such that $\cat  
M=\dim M$. Then $\Arn(M,\omega)\geq \Crit M$.
\endproclaim
\p The case $\dim M =2$ is well known, see [F89-1], [H88], so we  
assume that $\dim M \geq 4$. Consider any Hamiltonian  
symplectomorphism $\phi: M \to M$ and the corresponding data $\Phi$  
and $f: X \to M$ as in 3.5. Then, by 3.5, $\Fix \phi = \Rest \Phi$,  
and hence, by 3.4 and 2.5, $\Fix \phi \geq 1+\cat M$, and thus  
$\Arn(M,\omega)\geq 1+\cat M$. Furthermore, by a theorem of Takens  
[T68], $\Crit M\leq 1+\dim M$. Now,
$$
1+\cat M \leq \Crit M \leq 1+\dim M=1+\cat M,
$$
and thus $\Arn(M,\omega) \geq \Crit M$.
\qed
{\mytitle 3.7. Example {\rm $(\cat M > \cl M)$}.} Let $M$ be a  
four-dimensional aspherical symplectic manifold described in [MS95,  
Example 3.8]. It is easy to see that $H^1(M)=\Bbb Z^3$. Furthermore,  
$H^*(M)$ is torsion free, and so $a^2=0$ for every $a\in H^1(X)$.  
Hence, $\cl M=3$. However, $\cat M=4$  because $\cat V =\dim V$ for  
every closed aspherical manifold $V$, see [EG57]. Moreover, for  
every closed symplectic manifold $N$ we have $\cl (M \times N) <  
\cat (M \times N)$ because $\cat (M \times N)=\dim N +4$ accoring to  
[RO97].
\head {\bf \S 4. The invariant $r(M)$ and critical points} \endhead
\m Let $X$ be a $CW$-space and let $A,B$ be two $CW$-subspaces of  
$X$. Then for every spectrum $E$ we have the cap-product
$$
\cap : E_i(X, A\cup B) \otimes \Pi^j(X,A) \to E_{i-j}(X,B),
$$
see [Ad74], [Sw75]. Here $\Pi^*(-)$ denotes stable cohomotopy,  
i.e., $\Pi^*(-)$ is the cohomology theory represented by the sphere  
spectrum $S$.
\par In  particular, if $D=D^k$ is the $k$-dimensional disk then  
for every $CW$-pair $(X,A)$ we have the cup-product
$$
\cap: E_i(X \times D, X \times \partial D\cup A\times D) \otimes  
\Pi^k(X \times D, X \times \partial D) \to E_{i-k}(X \times D, A  
\times D).
$$
Let $a\in \Pi^k(D, \partial D)=\Bbb Z$ be a generator. We set  
$t=p^*a\in \Pi^k(X \times D, X \times \partial D)$ where $p: (X  
\times D, X \times \partial D) \to (D, \partial D)$ is the  
projection.
\proclaim{4.1. Lemma} For every $CW$-pair $(X,A)$ the homomorphism
 $$
\cap t: E_i(X \times D, X \times \partial D\cup A\times D) \to  
E_{i-k}(X \times D, A \times D)
$$
is an isomorphism.
\endproclaim
In fact, it is a relative Thom--Dold isomorphism.
\p If $A=\emptyset$ then $\cap t$ is the standard Thom--Dold  
isomorphism for the trivial $D^k$-bundle (or the suspension  
isomorphism, if you want), see e.g. [Sw75]. In other words, for  
$A=\emptyset$ the homomorphism in question has the form
$
\cap t: E_i(T \alpha) \to E_{i-k}(X \times D)
$
where $T \alpha$ is the Thom space of the trivial $D^k$-bundle  
$\alpha$. Furthermore, the homomorphism in question has the form
$$
E_i(T \alpha, T(\alpha|A)) \to E_{i-k}(X \times D, A \times D).
$$
Considering the commutative diagram
$$
\CD
\cdots \to @.  E_i(T(\alpha|A)) @>>> E_i(T \alpha) @>>> E_i(T  
\alpha, T(\alpha|A))@. \to \cdots\\
@. @VV\cap (t|A)V @VV\cap tV @VV\cap tV \\
\cdots \to @. E_{i-k}(A \times D) @>>> E_{i-k}(X \times D) @>>>  
E_{i-k}(X \times D, A \times D)@. \to  \cdots
\endCD
$$
with the exact rows, and using the Five Lemma, we conclude that the  
homomorphism in question is an isomorphism.
\qed
{\mytitle 4.2. Definition {\rm ([CZ83], [MO93])}.} Given a  
connected closed smooth manifold $M$, we define $\C G H_{p,q}(M)$ to  
be the set of all $C^2$-functions $g: M \times \Bbb R^{p+q} \to  
\Bbb R$ with the following properties:
\roster
\item There exist disks $D_+\i \Bbb R^p$ and $D_-\i \Bbb R^q$  
centered in origin such that $\inter (M \times D_+ \times D_-)$  
contains all critical points of $g$;
\item $\nabla g(x)$ points inward on $M \times \partial D_+ \times  
\inter D_-$ and outward on $M \times \inter D_+ \times \partial  
D_-$.
\endroster
{\mytitle 4.3. Definition {\rm ([CZ83], [MO93])}.} Given $g\in \C  
GH_{p,q}(M)$, consider the gradient flow $\dot x=\nabla g(x)$. Let  
$x\mul \Bbb R$ denote the solution of the flow through $x$. We  
choose $D_+$ and $D_-$ as in 4.2, set $B:=M \times D_+ \times D_-$  
and define $S_g=S_{g,B}:=\{x\in B|x\mul \Bbb R\i B\}$.
\proclaim{4.4. Theorem {\rm (cf [MO93, 4.1])}} For every function  
$g\in \C GH_{p,q}(M)$, there is a subpolyhedron $K$ of $\inter B$  
such that $S_g\i K$ and $\crit g\geq 1+\cat_BK$.
\endproclaim
\p We set $S=S_g$. Because of of 4.2, $S$ is a compact subset of  
$\inter B$. Furthermore, $S$ is an invariant set of the gradient  
flow $\dot x=\nabla g(x)$, and $S$ contains all critical points of  
$g$. Given $A\i S$, we define $\nu(A)=1+\cat_BA$. Clearly, $\nu$ is  
an index function on $S$. Thus, by 3.3, $\nu(S)\leq \crit g$. Now,  
let $V_1, \ldots V_{\nu(S)}$ be a covering of $S$ such that every  
$V_i$ is open and contractible in $B$. Choose any simplicial  
triangulation of $B$. Then, by the Lebesgue Lemma, there exists a  
simplicial subdivision of $B$ with the following property: every  
simplex $e$ with $e\cap S \neq \emptyset$ is contained in some  
$V_i$. Now, we set $K$ to be the union of all simplices $e$ with  
$e\cap S\neq  \emptyset$. Clearly, $1+\cat _BK \leq \nu(S)$, and  
thus $\crit g \geq 1+\cat _BK$. Finally, we can find $K\i \inter B$  
because of the collar theorem.
\qed
\m Let $r(M)$ be the invariant defined in 2.1.
\proclaim{4.5. Theorem} For every function $g\in \C GH_{p,q}(M)$,  
the number of critical points of $g$ is at least $1+r(M)$. In   
particular, $\crit g \geq 1+\cat M$ if $M$ is aspherical.
\endproclaim
\p Here we follow McCord--Oprea~[MO93]. However, unlike them, here  
we use certain extraordinary (co)homology instead of classical  
(co)homology.
\par Let $r:=r(M)$. We choose $K$ as in 4.4 and prove that  
$\cat_BK\geq r$. Consider the Puppe sequence
$$
P_r(M) @>p_r>>  M @>j_r>>  C_r(M).
$$
Let $e: M^+ \to C_r(M)$ be a map such that $e|M=j_r$ and $e$ maps  
the added point to the base point of $C_r(M)$. Let $h: C_r(M) \to C$  
be a pointed homotopy equivalence such that $C$ is a $CW$-complex.  
We set $E= \Sigma^{\infty}C$ and let $u_r\in E^0(M)$ be the stable  
homotopy class of the map $he: M^+ \to C$. Then $u_r \neq 0$ since  
$j_r$ is stably essential. \par
We define
$$
f: K \i B =M \times \Bbb R^{p+q}@>\text{projection}>> M.
$$
\proclaim{4.6. Lemma} If $f^*u_r \neq 0$ then $\cat_BK\geq r$.
\endproclaim
\p  Since $(p_r, j_r)$ is a Puppe sequence, $p_r^*u_r=0$. Hence,  
the map $f$ can't be lifted to $P_r(M)$, and therefore the inclusion  
$K\i B$ can't be lifted to $P_r(B)$. So, $\cat _BK\geq r$. The  
lemma is proved.
\m We continue the proof of the theorem. Let $j: K \i B$ be the  
inclusion. By 4.6, it suffices to prove that $j^*:E^*(B) \to E^*(K)$  
is a monomorphism. Notice that if $Y$ is a $CW$-subspace of $\Bbb  
R^N$ then there is a duality isomorphism
$$
E^0(Y) \cong E_{-N}(\Bbb R^N, \Bbb R^N\setminus Y):=E_{-N}(\Bbb  
R^N\cup C(\Bbb R^N\setminus Y))
$$
see e.g. [DP84]. So, it suffices to prove that the dual homomorphism 
$$
D(j^*): E_*(\Bbb R^N, \Bbb R^N\setminus B) \to E_*(\Bbb R^N, \Bbb  
R^N \setminus K)
$$
is monic for a certain (good) embedding $B \to \Bbb R^N$.
\par  We have the following commutative diagram:
$$
\CD
E_*(\Bbb R^N, \Bbb R^N\setminus B) @= E_*(\Bbb R^N, \Bbb  
R^N\setminus B)\\
@VhV\cong V @VD(j^*)VV\\
E_*(\Bbb R^N, \Bbb R^N\setminus \inter B) @>>> E_*(\Bbb R^N, \Bbb  
R^N\setminus K)\\
@AeA\cong A @Ae'A\cong A\\
E_*(B, \partial B) @>a_*>> E_*(B, B\setminus K)
\endCD
$$
where all the homomorphisms except $D(j^*)$ are induced by the  
inclusions. Here $h$ is an isomorphism since the inclusion $\inter B  
\to B$ is a  homotopy equivalence (the space $B\setminus$  
\{collar\} is a deformation retract of $\inter B$). Furthermore, $e$  
is an isomorphism since $(B, \partial B)$ and $(\Bbb R^N, \Bbb  
R^N\setminus \inter B)$ are cofibered pairs, while $e'$ is an  
isomorphism by Lemma 3.4 from [DP84]. So, $D(j^*)$ is monic if $a_*$  
is. Since $B\setminus K\i B\setminus S$, it suffices to prove that   
$E_*(B, \partial B) \to  E_*(B, B\setminus S)$ is a monomorphism.
\par
Let $B_+=M \times \partial D_+ \times D_-$, and let $B_-=M \times  
D_+ \times \partial D_-$. Furthermore, let $A_+:=\{x\in B\bigm|x\mul  
\Bbb R_-\in B\}$ and let $A_-:=\{x\in B\bigm|x\mul \Bbb R_+\in  
B\}$. Then $B_+ \cap A_-=\emptyset=B_-\cap A_+$, and so there are  
the inclusions $i_+: (B,B_+) \to (B,B\setminus A_-)$ and $i_-:  
(B,B_-) \to (B,B\setminus A_+)$. It turns out to be that both $i_+$  
and $i_-$ are homotopy equivalences, [CZ83, Lemma 3].
\par Let $t\in \Pi^m(B, B_-)$ be the class as in 4.1, and let  
$t':=((i_-)^*)^{-1}(t)$. Since $S=A_+ \cap A_-$, we have the  
commutative diagram
$$
\CD
E_i(B, \partial B) @>>> E_i(B, B\setminus S)\\
@V \cong V\cap tV @VV\cap t'V\\
E_{i-q}(B,B_+) @> \cong >> E_{i-q}(B,B\setminus A_-)
\endCD
$$
where the left map is an isomorphism by 4.1 and the bottom map is  
the isomorphism $(i_+)_*$. (Generally, $(B,B\setminus S)$ is not a  
$CW$-pair, but nevertheless in our case the map $\cap t'$ is  
defined, see [DP84, 3.5].)
Thus, the top homomorphism is injective.
\par Finally, if $M$ is aspherical then $\cat M =\dim M$, EG57],  
and so $r(M)=\cat M$ by  2.4.
\qed
\head{References}\endhead
\hyphenation{To-pol-ogy}
\halign{{\bf #\ }\hfil & \vtop{\parindent0pt
\hsize=31.1em
\hangindent0em\strut#\strut}\cr
[Ad74] & J.F. Adams: \it Stable Homotopy and Generalised  
Cohomology. \rm The Univ. of Chicago Press, Chicago 1974\cr
[Ar89] & V.I. Arnold: \it Mathematical Methods in Classical  
Mechanics. \rm Springer, Berlin Heidelberg New York 1989\cr
[CP86] & M. Clapp, D. Puppe: \it Invariants of  
Lusternik--Schnirelmann type and the topology of critical sets. \rm  
Trans. Amer. Math. Soc. {\bf 298}, 2 (1986), 603--620 \cr
[C76] & C. Conley: \it Isolated invariant sets and the Morse index.
\rm CBMS Regional Conf. Ser. in Math. {\bf 38}, Amer. Math. Soc.,  
Providence, R. I., 1976 \cr
[CZ83] & C. Conley, E. Zehnder: \it The Birkhoff--Lewis Fixed Point  
Theorem and a Conjecture of V.I. Arnold, \rm Invent. Math. {\bf 73}  
(1983), 33-49\cr
[Co98] & O. Cornea: \it Some properties of the relative  
Lusternik--Schnirelmann category. \rm Proc of the homotopy year of  
the Fields Institute, to appear\cr
[DP84] & A. Dold, D. Puppe: \it Duality, trace and transfer. \rm  
Proceedings of the Steklov Inst. of Math. {\bf 4} (1984), 85--103  
\cr
[EG57] & S. Eilenberg and T. Ganea: \it On the  
Lusternik--Schnirelmann category of abstract groups, \rm Ann. of  
Math. {\bf 65} (1957), 517--518\cr
[F89-1] & A. Floer: \it Cuplength estimates on Lagrangian  
intersections, \rm Commun. Pure Appl. Math {\bf 42} (1989),  
335--356\cr
[F89-2] & A. Floer: \it Symplectic fixed points and holomorphic  
spheres, \rm Commun. Math. Phys. {\bf 120} (1989), 575--611\cr
[Fox41] & R. Fox: \it On the Lusternik--Schnirelmann category\rm ,  
Ann. of Math. {\bf 42} (1941), 333-370\cr
[FP90] & R. Fritsch, R.A. Piccinini: {\it Cellular structures in  
topology.} Cambridge Univ. Press, Cambridge, 1990\cr
[H88] & H. Hofer: \it Lusternik--Schnirelmann theory for Lagrangian  
intersections. \rm Annales de l'inst. Henri Poincar\'e-- analyse  
nonlineare, {\bf  5} (1988), 465--499\cr
[HZ] & H. Hofer, E. Zehnder: Symplectic Invariants and
Hamiltonian Dynamics, Birkh\"auser, Basel, 1994\cr
[Hu61] & P. Huber: Homotopical cohomology and \v Cech cohomology.  
Math. Annalen {\bf 144} (1961), 73--76\cr
[J78] & I.M. James: \it On category, in the sense of  
Lusternik--Schnirelmann. \rm Topology, {\bf 17} (1978), 331--349\cr
[LO96] & H.V Le, K. Ono: {\it Cup-length estimate for symplectic  
fixed points.} in: (Contact and Symplectic Geometry (ed. by C.  
Thomas), London Math. Soc. Lecture Notes. Ser., Cambridge Univ.  
Press, 1996\cr
[MO93] & C. McCord, J. Oprea: \it Rational Lusternik--Schnirelmann  
category and the Arnol'd conjecture for nilmanifolds. \rm Topology,  
{\bf 32}, 4 (1993), 701--717\cr
[MS95] & D. McDuff, D. Salamon: \it Introduction to Symplectic  
Topology, \rm Clarendon Press, Oxford 1995\cr
[M59] & J. Milnor: {\it On spaces having the homotopy type of a  
$CW$-complex}. Trans. Amer. Math. Soc. {\bf 90} (1959), 272--280\cr
[Mu66] &J. Munkres \it  Elementary Differential Topology, \rm Ann.  
of Math. Studies {\bf 54}, Princeton Univ. Press, Princeton, New  
Jersey 1966 \cr
[R98] & Yu. B. Rudyak: \it On category weight and its applications,  
\rm to appear in Topology, (1998)\cr
[RO97] & Yu. B. Rudyak, J. Oprea: \it On the  
Lusternik--Schnirelmann Category of Symplectic
Manifolds and the Arnold Conjecture, \rm Preprint, available as  
dg-ga/9708007 \cr
[S85] & J.-C. Sikorav: {\it Points fixes d'une application  
symplectique.} Jour. Diff. Geom., {\bf 25} (1985), 49--79 \cr
[Sp66] & E.H. Spanier: \it Algebraic Topology. McGraw-Hill, New  
York 1966.\cr
[St68] & J. Stallings: \it Lectures on Polyhedral Topology. \rm  
Tata Institute of Fundamental Research, Bombay 1968\cr
[Sv66] & A. Svarc: \it The genus of a fiber space. \rm Amer. Math.  
Soc. Translations {\bf 55} (1966), 49--140\cr
[Sw75] & R.W. Switzer: \it Algebraic Topology -- Homotopy and  
Homology. \rm Springer, Berlin Heidelberg New York 1975\cr
[T68] & F. Takens: \it The minimal number of crirical points of a  
function on a compact manifold and the Lusternik--Schnirelmann  
category, \rm Invent. Math. {\bf 6} (1968), 197--244 \cr
\cr [W70] & C.T.C. Wall: \it Surgery on compact manifolds. \rm  
Academic Press, New York \& London 1970\cr}
\enddocument